# NOTES ON THE TWO-DIMENSIONAL FRACTIONAL BROWNIAN MOTION

BY FABRICE BAUDOIN AND DAVID NUALART[1]

*Toulouse University of Technology and Universitat de Barcelona*

We study the two-dimensional fractional Brownian motion with Hurst parameter $H > \frac{1}{2}$. In particular, we show, using stochastic calculus, that this process admits a *skew-product* decomposition and deduce from this representation some asymptotic properties of the motion.

**1. Introduction.** We consider a complex fractional Brownian motion started at $z_0 = x^1 + ix^2 \in \mathbb{C}^*$,

$$B_t = B_t^1 + iB_t^2, \qquad t \geq 0,$$

where $B^1$ and $B^2$ are two independent fractional Brownian motions with the same parameter $H \in (0,1)$. That is, for $i = 1, 2$, $B^i = \{B_t^i, t \geq 0\}$ is a Gaussian process with mean $\mathbb{E}(B_t^i) = x^i$ and covariance

$$R(t,s) = \tfrac{1}{2}(s^{2H} + t^{2H} - |t-s|^{2H}).$$

Notice that, in the case $H = \frac{1}{2}$, $B$ is a planar Brownian motion. In what follows, we shall always assume that $H > \frac{1}{2}$. In this case, it is known that $B$ is a transient process (see [20]).

The aim of this paper is to use the stochastic calculus for the fractional Brownian motion which has been recently developed by several authors (see, e.g., [1, 2, 4, 11]) in order to obtain geometric properties of this motion as it has been made in the case of the planar Brownian motion (see, e.g., [6, 8, 12, 13, 17, 18]). The key point of our study is an analogue of the celebrated skew-product decomposition of the two-dimensional Brownian motion. Precisely, we show that we can write

$$B_t = z_0 \exp\left(\int_0^t \frac{dB_s}{B_s}\right), \qquad t \geq 0,$$

Received March 2003; revised October 2004.
[1]Supported by DGES Grant BFM2000-0598.

where the integral of the right-hand side can be understood not only pathwise, but also as a local divergence. We shall see that it is then possible to apply the ergodic theorem to study the asymptotics of the integral

$$\int_0^t \frac{dB_s}{B_s}, \tag{1.1}$$

whose real part (resp. the imaginary part) explains the radial part of $B$ (resp. the angular part of $B$).

In that respect, the study of the two-dimensional fractional Brownian motion with Hurst parameter $H > \frac{1}{2}$ could seem simpler than the study of the planar Brownian motion, for which it is not possible to apply directly the ergodic theorem. Nevertheless, for the fractional Brownian motion, we shall see that the study of the windings is much more difficult because the integral (1.1) is not a time-changed fractional Brownian motion.

**2. Itô's formula for holomorphic functions.** Before we study the skew-product decomposition, we need to develop stochastic calculus for the two-dimensional fractional Brownian motion. In particular, we shall need a two-dimensional Itô's formula.

Let $\mathcal{H}$ be the Hilbert space defined as the closure of the set of step functions on $\mathbb{R}_+$ with respect to the scalar product

$$\langle \mathbf{1}_{[0,t]}, \mathbf{1}_{[0,s]} \rangle_\mathcal{H} = R(t,s).$$

For $i = 1, 2$, the mapping $\mathbf{1}_{[0,t]} \longrightarrow B_t^i - x^i$ can be extended to an isometry between $\mathcal{H}$ and the first chaos $H_1^i$ associated with $B^i$. We denote this isometry by $\varphi \to B^i(\varphi)$. If $\varphi \in \mathcal{H}^{\otimes 2}$, we set $B(\varphi) = (B^1(\varphi^1), B^2(\varphi^2))$.

Let $\mathcal{S}$ be the set of smooth and cylindrical random variables of the form

$$F = f(B(\varphi_1), \ldots, B(\varphi_n)), \tag{2.1}$$

where $n \geq 1$, $f \in C_b^\infty(\mathbb{R}^{2n})$ ($f$ and all its partial derivatives are bounded), and $\varphi_i \in \mathcal{H}^{\otimes 2}$. The derivative operator $D$ of a smooth and cylindrical random variable $F$ of the form (2.1) is defined as the $\mathcal{H}^{\otimes 2}$-valued random variable

$$DF = \sum_{j=1}^n \left( \frac{\partial f}{\partial x_j}(Z) \varphi_j^1, \frac{\partial f}{\partial y_j}(Z) \varphi_j^2 \right),$$

where $Z = (B(\varphi_1), \ldots, B(\varphi_n))$, and we make use of the notation $f(z_1, \ldots, z_n)$ and $z_j = (x_j, y_j)$, for $j = 1, \ldots, n$. The derivative operator $D$ is then a closable operator from $L^p(\Omega)$ into $L^p(\Omega; \mathcal{H}^{\otimes 2})$ for any $p \geq 1$. For any $p \geq 1$, the Sobolev space $\mathbb{D}^{1,p}$ is the closure of $\mathcal{S}$ with respect to the norm

$$\|F\|_{1,p}^p = \mathbb{E}|F|^p + \mathbb{E}\|DF\|_{\mathcal{H}^{\otimes 2}}^p.$$

In a similar way, given a Hilbert space $V$, we denote by $\mathbb{D}^{1,p}(V)$ the corresponding Sobolev space of $V$-valued random variables.

The divergence operator $\delta$ is the adjoint of the derivative operator, defined by means of the duality relationship

$$\mathbb{E}(F\delta(u)) = \mathbb{E}\langle DF, u\rangle_{\mathcal{H}^{\otimes 2}}, \tag{2.2}$$

where $u$ is a random variable in $L^2(\Omega; \mathcal{H}^{\otimes 2})$. That is, we say that $u$ belongs to the domain of the operator $\delta$, denoted by $\mathrm{Dom}\,\delta$, if the mapping $F \mapsto E\langle DF, u\rangle_{\mathcal{H}^{\otimes 2}}$ is continuous in $L^2(\Omega)$. If $u \in \mathrm{Dom}\,\delta$ is a two-dimensional stochastic process, we will make use of the notation

$$\delta(u) = \int_0^\infty u_s^1\,\delta B_s^1 + \int_0^\infty u_s^2\,\delta B_s^2.$$

The space $\mathbb{D}^{1,2}(\mathcal{H}^{\otimes 2})$ is included in $\mathrm{Dom}\,\delta$. For an element $u$ in $\mathbb{D}^{1,2}(\mathcal{H}^{\otimes 2})$, we have

$$\mathbb{E}(\delta(u))^2 = \mathbb{E}\|u\|_{\mathcal{H}^{\otimes 2}}^2 + \mathbb{E}\langle Du, (Du)^*\rangle_{\mathcal{H}^{\otimes 4}}. \tag{2.3}$$

We recall that $\mathcal{H}$ contains the space of real-valued functions $\varphi$ on $\mathbb{R}_+$ such that

$$\int_0^\infty \int_0^\infty |\varphi(t)||\varphi(s)||t-s|^{2H-2}\,dt\,ds < \infty.$$

The following result has been proved in [2].

PROPOSITION 1.   *Let $u = (u_t)_{t\geq 0}$ be a two-dimensional stochastic process in the space $\mathbb{D}^{1,2}(\mathcal{H}^{\otimes 2})$ such that, for some $T > 0$,*

$$\mathbb{E}\left(\int_0^T \int_0^T |u_t||u_s||t-s|^{2H-2}\,dt\,ds\right) < \infty,$$

$$\mathbb{E}\int_{[0,T]^4} |D_t u_\theta||D_s u_\eta||t-s|^{2H-2}|\theta-\eta|^{2H-2}\,dt\,du\,d\theta\,d\eta < \infty$$

*and*

$$\int_0^T \int_0^T |D_s u_t||t-s|^{2H-2}\,ds\,dt < \infty,$$

*a.s. Then the symmetric integral defined as the limit in probability*

$$\int_0^T u_t\,dB_t = \lim_{\varepsilon\downarrow 0}(2\varepsilon)^{-1}\int_0^T u_s(B_{(s+\varepsilon)\wedge T} - B_{(s-\varepsilon)\vee 0})\,ds$$

*exists and we have*

$$\int_0^T u_t\,dB_t = \delta(u) + \alpha_H \int_0^T \int_0^T (D_s^1 u_t^1 + D_s^2 u_t^2)|t-s|^{2H-2}\,ds\,dt, \tag{2.4}$$

*where $\alpha_H = H(2H-1)$.*

We denote by $\mathbb{D}^{1,2}_{\mathrm{loc}}(\mathcal{H}^{\otimes 2})$ the set of $\mathcal{H}^{\otimes 2}$-valued random variables $u$ such that there exists a sequence $\{(\Omega^n, u^n)\}, n \geq 1\} \subset \mathcal{F} \times \mathbb{D}^{1,2}(\mathcal{H}^{\otimes 2})$ such that $\Omega^n \uparrow \Omega$ a.s. and $u = u^n$, a.e. on $[0,T] \times \Omega_n$. We then say that $\{(\Omega^n, u^n)\}$ localizes $u$ in $\mathbb{D}^{1,2}(\mathcal{H}^{\otimes 2})$. If $u \in \mathbb{D}^{1,2}_{\mathrm{loc}}(\mathcal{H}^{\otimes 2})$, by the local property of the divergence operator, we can define without ambiguity $\delta(u)$ by setting

$$\delta(u)|_{\Omega^n} = \delta(u^n)|_{\Omega^n}$$

for each $n \geq 1$, where $\{(\Omega^n, u^n)\}$ is a localizing sequence for $u$ in $\mathbb{D}^{1,2}(\mathcal{H}^{\otimes 2})$.

If $f$ is a real-valued function on $\mathbb{R}^2$ twice continuously differentiable, then for each $t > 0$, $\nabla f(B.) \mathbf{1}_{(0,t)}$ belongs to $\mathbb{D}^{1,2}_{\mathrm{loc}}(\mathcal{H}^{\otimes 2})$ and the following version of Itô's formula holds (see [2]):

$$
\begin{aligned}
(2.5) \quad f(B_t) &= f(x^1, x^2) + \int_0^t \frac{\partial f}{\partial x}(B_s) \, \delta B_s^1 + \int_0^t \frac{\partial f}{\partial y}(B_s) \, \delta B_s^2 \\
&\quad + H \int_0^t \left( \frac{\partial^2 f}{\partial x^2}(B_s) + \frac{\partial^2 f}{\partial y^2}(B_s) \right) s^{2H-1} \, ds.
\end{aligned}
$$

On the other hand, we will denote by $\int_0^t u_s \, dB_s$ the Riemann–Stieltjes integral of a stochastic process $u_t$ whose trajectories are $\gamma$-Hölder continuous, provided $\gamma > 1 - H$ (see [21]).

If $f : \mathbb{C} \to \mathbb{C}$ is a twice continuously differentiable function in the variables $(x, y) \to f(x + iy)$, formula (2.5) leads to

$$
\begin{aligned}
(2.6) \quad f(B_t) &= f(z_0) + \int_0^t \frac{\partial f}{\partial z}(B_s) \, \delta B_s^1 \\
&\quad + \int_0^t \frac{\partial f}{\partial \overline{z}}(B_s) \, \delta \overline{B}_s^2 + H \int_0^t \Delta f(B_s) s^{2H-1} \, ds,
\end{aligned}
$$

where $\frac{\partial f}{\partial z} = \frac{1}{2}(\frac{\partial f}{\partial x} - i\frac{\partial f}{\partial y})$ and $\frac{\partial f}{\partial \overline{z}} = \frac{1}{2}(\frac{\partial f}{\partial x} + i\frac{\partial f}{\partial y})$.

PROPOSITION 2. *Let $f : \mathbb{C} \to \mathbb{C}$ be a holomorphic function, then*

$$(2.7) \qquad f(B_t) = f(z_0) + \int_0^t f'(B_s) \, \delta B_s$$

$$(2.8) \qquad\qquad = f(z_0) + \int_0^t f'(B_s) \, dB_s,$$

*where $f'(z) = \frac{\partial f}{\partial z}$.*

PROOF. Equation (2.7) follows from (2.5) and the fact that $\frac{\partial f}{\partial \overline{z}} = 0$ and $\Delta f = 0$. Equation (2.8) follows from the change of variables formula for integrals with respect to Hölder-continuous functions. $\square$

REMARK 3. Unlike two-dimensional Brownian motion, two-dimensional fractional Brownian motion is not conformal invariant. It means that, in general, $\int_0^t f'(B_s) \delta B_s$ is not a time-changed fractional Brownian motion. Nevertheless, let us recall that, in [3], it is shown that, for a deterministic function, $g$ which satisfies some regularity assumptions and for a one-dimensional fractional Brownian motion $B^H$ with Hurst parameter $H > \frac{1}{2}$, we can write a decomposition

$$\int_0^t g(s)\,dB_s^H = \beta^H_{\int_0^t g(s)^{1/H}\,ds} + A_t, \qquad t \geq 0,$$

where $A_t$ is a bounded variation process and $\beta^H$ a fractional Brownian motion with Hurst parameter $H$. Hence, a natural question is to ask if such a decomposition exists (with the suitable time-change) for the two-dimensional process $f(B_t)$ or not, but so far there is no answer yet to this question.

As an immediate corollary of the previous proposition, we also deduce the following:

COROLLARY 4. Let $\alpha = f(z)\,dz$ be a holomorphic one-form on $\mathbb{C}$, then

$$\int_{B(0,t)} \alpha = \int_0^t f'(B_s) \delta B_s = \int_0^t f'(B_s)\,dB_s,$$

where $\int_{B(0,t)} \alpha$ denotes the integral of $\alpha$ along the trajectories of $B$.

REMARK 5. We recall that we can define the integral of a holomorphic one-form on any continuous path.

REMARK 6. The previous equality is interesting because it relates objects of different natures.

From now on, since the different integrals coincide for holomorphic functions, we shall use indifferently the notation $\delta$ or $d$.

As a consequence of Itô's formula, we can deduce the following.

PROPOSITION 7. Let $f : \mathbb{C} \to \mathbb{C}$ be a holomorphic function, then

$$\mathbb{E}\left(\left(\operatorname{Im}\int_0^t f(B_s)\,dB_s\right)^2\right) = \mathbb{E}\left(\left(\operatorname{Re}\int_0^t f(B_s)\,dB_s\right)^2\right)$$

and

$$\mathbb{E}\left(\left(\operatorname{Im}\int_0^t f(B_s)\,dB_s\right)\left(\operatorname{Re}\int_0^t f(B_s)\,dB_s\right)\right) = 0.$$

PROOF. In fact, this is a consequence of Itô's formula. Indeed, let us consider a holomorphic function $F$ such that $F = f'$. From Itô's formula, we have to show that

$$\mathbb{E}((\mathrm{Im}(F(B_t) - F(z_0)))^2) = \mathbb{E}((\mathrm{Re}(F(B_t) - F(z_0)))^2)$$

and

$$\mathbb{E}((\mathrm{Im}(F(B_t) - F(z_0)))(\mathrm{Re}(F(B_t) - F(z_0)))) = 0.$$

Now, notice that after a suitable scaling, these equalities only involve the centered reduced normal law, and are, hence, easily checked. □

The following proposition provides the skew-product decomposition for the fBm.

PROPOSITION 8 ("Skew-product decomposition for the fBm"). *We have*

$$B_t = z_0 \exp\left(\int_0^t \frac{\delta B_s}{B_s}\right) = z_0 \exp\left(\int_0^t \frac{dB_s}{B_s}\right), \qquad t \geq 0.$$

PROOF. If we apply Itô's formula to

$$B_t \exp\left(-\int_0^t \frac{dB_s}{B_s}\right),$$

we see that it is constant, and so equal to $z_0$. □

Notice the important fact that the divergence integral $\int_0^t \frac{\delta B_s}{B_s}$ is local. Indeed, if not, we would have

$$\mathrm{Re}\left(\mathbb{E}\left(\int_0^t \frac{\delta B_s}{B_s}\right)\right) = \mathbb{E}\left(\int_0^t \frac{B_s^1 \delta B_s^1 + B_s^2 \delta B_s^2}{|B_s|^2}\right) = 0,$$

and from Itô's formula,

$$\mathbb{E}(\ln|B_t|) = \mathbb{E}(\ln|z_0|),$$

which is clearly absurd.

**3. Polar decomposition.** First of all, we study the polarity of points for the complex fractional Brownian motion $B$ issued from $z_0 \neq 0$. For this, it is enough to adapt the original proof by [6] of the polarity of points for the planar Brownian motion.

PROPOSITION 9. *If $H \geq \frac{1}{2}$, then the points are polar for $B$, that is,* $\forall x \in \mathbb{C}$,

$$\mathbb{P}(\exists t > 0, B_t = x) = 0.$$

PROOF. We proceed in several steps. We assume $H > \frac{1}{2}$ because the case $H = \frac{1}{2}$, which corresponds to the case of the Brownian motion, is well known.

In [19] it has proved that a.s. the Hausdorff measure $\mu_\varphi$ of the fractional Brownian curve is finite, where $\varphi(\varepsilon) = \varepsilon^{1/H} \log\log(1/\varepsilon)$. This clearly implies that the Lebesgue measure of the fractional Brownian curve is zero. For the sake of completeness, we provide a simple and direct proof of this result. Let us denote by $\mu$ the Lebesgue measure on the plane and note that

$$\mathbb{E}(\mu\{B_t(\omega), 0 \leq t \leq 1\}) \leq \pi \mathbb{E}\left(\sup_{0 \leq t \leq 1} |B_t|^2\right).$$

Hence,

$$\mathbb{E}(\mu\{B_t(\omega), 0 \leq t \leq 1\}) < +\infty.$$

Now, from the scaling property of the fBm, for $T \in \mathbb{N}^*$,

$$\mathbb{E}(\mu\{B_t(\omega), 0 \leq t \leq T\}) = T^{2H} \mathbb{E}(\mu\{B_t(\omega), 0 \leq t \leq 1\}).$$

But, since the increments are stationary,

$$\mathbb{E}(\mu\{B_t(\omega), 0 \leq t \leq T\}) \leq \sum_{i=1}^{T} \mathbb{E}(\mu\{B_t(\omega), i-1 \leq t \leq i\})$$
$$\leq T \mathbb{E}(\mu\{B_t(\omega), 0 \leq t \leq 1\}).$$

We deduce that

$$T^{2H} \mathbb{E}(\mu\{B_t(\omega), 0 \leq t \leq 1\}) \leq T \mathbb{E}(\mu\{B_t(\omega), 0 \leq t \leq 1\}),$$

which implies

$$\mathbb{E}(\mu\{B_t(\omega), t \geq 0\}) = 0.$$

Now, if there exists at least one $x \in \mathbb{C} \setminus \{z_0\}$ such that

$$\mathbb{P}(\exists t > 0, B_t = x) > 0,$$

then we claim that, for all $z \in \mathbb{C} \setminus \{z_0\}$,

$$\mathbb{P}(\exists t > 0, B_t = z) > 0.$$

This follows from the rotational invariance and the scaling property of the fBm. Since

$$\mathbb{E}(\mu\{B_t(\omega), t \geq 0\}) = 0,$$

we have, for all $z \in \mathbb{C} \setminus \{z_0\}$,

$$\mathbb{P}(\exists t > 0, B_t = z) = 0.$$

To conclude the proof, we have now to show that
$$\mathbb{P}(\exists t > 0, B_t = z_0) = 0.$$
For this, notice that for, $h > 0$, the process $(B_{t+h} - B_h)_{t \geq 0}$ is still a fractional Brownian motion (this is a Gaussian process with the same covariance function as $B$). Hence, if we had
$$\mathbb{P}(\exists t > 0, B_t = z_0) > 0,$$
then, we would have, for $h > 0$,
$$\mathbb{P}(\exists t > 0, B_{t+h} - B_h = z_0) > 0,$$
which is absurd because $z_0 \neq 0$. $\square$

REMARK 10. It is known (see, e.g., [20]) that, for $H \in (0, \frac{1}{2})$, the sample paths of $B$ have interior points, so that $\mathbb{E}(\mu\{B_t(\omega), 0 \leq t \leq 1\}) > 0$. In particular, the points are not polar. Notice also that, for $H \in (0, \frac{1}{2}]$, $B$ is recurrent, whereas it is transient for $H > \frac{1}{2}$.

Assume now that $z_0 = 0$. Since $0$ is polar for the complex fBm, we can consider for $t > 0$ the radial part
$$\rho_t = |B_t|$$
and the angular part
$$\Theta_t = \frac{B_t}{|B_t|}.$$

PROPOSITION 11. *For each $t > 0$, $\Theta_t$ is independent of the process $(\rho_s)_{s \geq 0}$ and uniformly distributed on the unit circle.*

PROOF. Let $X$ be an $\mathcal{R}_\infty$-measurable random variable, where $\mathcal{R}$ is the natural filtration of $(\rho_s)_{s \geq 0}$ and let $f$ be a Borel function on the circle. We know that the law of $B$ is rotational invariant, and this implies that, for all $\theta \in [0, 2\pi]$,
$$\mathbb{E}(Xf(\Theta_t)) = \mathbb{E}(Xf(e^{i\theta}\Theta_t)).$$
Integrating this relation yields
$$\mathbb{E}(Xf(\Theta_t)) = \mathbb{E}\left(X \frac{1}{2\pi} \int_0^{2\pi} f(e^{i\theta}\Theta_t) \, d\theta\right).$$
Since the Lebesgue measure on the circle is invariant by translation, we deduce that
$$\frac{1}{2\pi} \int_0^{2\pi} f(e^{i\theta}\Theta_t) \, d\theta$$

is constant and does not depend on $\Theta_t$. This implies that

$$\mathbb{E}(Xf(\Theta_t)) = \mathbb{E}(X)\frac{1}{2\pi}\int_0^{2\pi} f(e^{i\theta})\,d\theta,$$

which gives the expected result. $\square$

**4. Asymptotics of the skew-product and limit theorems.** We come back to the case $z_0 \neq 0$ and study the asymptotics of $\int_0^t \frac{\delta B_s}{B_s} = \int_0^t \frac{dB_s}{B_s}$ as $t$ tends to infinity. First we present some preliminary results.

LEMMA 12. *For any real number $k > 1$, we have*

$$\mathbb{E}\left(\left|\int_1^k \frac{dB_s}{B_s - z_0}\right|\right) < +\infty.$$

PROOF. In order to show this property, we express the path-wise Riemann–Stieltjes integral in terms of fractional derivatives, using the fractional integration by parts formula (see [22]):

$$\int_1^k \frac{dB_s}{B_s - z_0} = \int_1^k D_{1+}^\alpha\left(\frac{1}{B_\cdot - z_0}\right)(s) D_{k-}^{1-\alpha} B_{k-}(s)\,ds,$$

where $1 - H < \alpha < \frac{1}{2}$, and the left and right fractional derivatives are given by

$$D_{k-}^{1-\alpha} B_{k-}(s) = \frac{-1}{\Gamma(\alpha)}\left(\frac{B_k - B_s}{(k-s)^{1-\alpha}} + \alpha\int_s^k \frac{B_y - B_s}{(y-s)^{2-\alpha}}\,dy\right)$$

and

$$D_{1+}^\alpha\left(\frac{1}{B_\cdot - z_0}\right)(s)$$
$$= \frac{1}{\Gamma(1-\alpha)}\left(\frac{(B_s - z_0)^{-1}}{(s-1)^\alpha} + \alpha\int_1^s \frac{1/(B_s - z_0) - 1/(B_y - z_0)}{(s-y)^{\alpha+1}}\,dy\right).$$

Clearly, for all $p \geq 1$,

(4.1) $$\mathbb{E}\left(\sup_{1\leq s\leq k} |D_{k-}^{1-\alpha} B_{k-}(s)|^p\right) < \infty.$$

On the other hand, for $1 < p < 2$,

(4.2) $$\mathbb{E}\left(\int_1^k \frac{ds}{|B_s - z_0|^p (s-1)^{\alpha p}}\right) < \infty$$

and

$$\mathbb{E}\bigg(\int_1^k \bigg(\int_1^s \frac{|1/(B_s - z_0) - 1/(B_y - z_0)|}{(s-y)^{\alpha+1}} dy\bigg)^p ds\bigg)$$

$$\leq \mathbb{E}\bigg(\int_1^k \bigg(\int_1^s \frac{|B_y - B_s|}{(s-y)^{\alpha+1}|B_s - z_0||B_y - z_0|} dy\bigg)^p ds\bigg)$$

(4.3)

$$\leq k^p \mathbb{E}\bigg[G^p \int_1^k \int_1^s \frac{1}{(s-y)^{p(\alpha+1-H+\varepsilon)}|B_s - z_0|^p |B_y - z_0|^p} dy\,ds\bigg]$$

$$\leq Ck^p (\mathbb{E}(G^{pq'}))^{1/q'} \sup_{1\leq s\leq k} \bigg(\mathbb{E}\bigg(\frac{1}{|B_s - z_0|^{2pq}}\bigg)\bigg)^{2/q} < \infty,$$

where

$$G = \sup_{1\leq y<s\leq k} \frac{|B_y - B_s|}{(s-y)^{H-\varepsilon}},$$

provided $\frac{1}{q} + \frac{1}{q'} = 1$, $2pq < 2$, and $p(\alpha + 1 - H + \varepsilon) < 1$. The estimates (4.1), (4.2) and (4.3) imply the desired result. □

REMARK 13. The previous result is not true for $H = \frac{1}{2}$, that is, in the case of the planar Brownian motion.

Now, we prove the ergodicity of the scaling for the fractional Brownian motion.

LEMMA 14. *For $k \in \mathbb{R}_+^* \setminus \{1\}$, the transformation on the path space $T_k : \omega \to \frac{\omega(k\cdot)}{k^H}$ which preserves the law of a one-dimensional fBm issued from 0 is ergodic.*

PROOF. Let $\beta$ be a one-dimensional fractional Brownian motion issued from 0. It suffices to show that $T_k$ is mixing. For this, we note that, for $s,t > 0$ and $n \in \mathbb{N}^*$,

$$\mathbb{E}\bigg(\beta_s \frac{\beta_{k^n t}}{k^{nH}}\bigg) = \frac{1}{2k^{nH}}(s^{2H} + k^{2nH}t^{2H} - |s - k^n t|^{2H})$$

so that

$$\lim_{n\to+\infty} \mathbb{E}\bigg(\beta_s \frac{\beta_{k^n t}}{k^{nH}}\bigg) = 0,$$

which implies the mixing property of $T_k$. □

Before we turn to the skew-product decomposition, we state a last preliminary result on the $\frac{1}{H}$-variation of the integral $\int_0^t \frac{dB_s}{B_s}$. We refer to [5] for similar results on the $\frac{1}{H}$-variation of divergence integrals with application to fractional Bessel processes.

PROPOSITION 15. *Set $t_i = \frac{it}{n}$, for $i = 0, 1, \ldots, n$. Then*

$$\sum_{i=0}^{n-1} \left| \int_{t_i}^{t_{i+1}} \frac{dB_s}{B_s} \right|^{1/H} \to c \int_0^t \frac{ds}{|B_s|^{1/H}},$$

*almost surely, as $n$ tends to infinity, where $c = E(|B_1|^{1/H})$.*

PROOF. For any complex-valued process $u_t$ whose paths are locally $\gamma$-Hölder continuous, with $\gamma > 1 - H$, we set

$$V_n(u) = \sum_{i=0}^{n-1} \left| \int_{t_i}^{t_{i+1}} u_s \, dB_s \right|^{1/H}.$$

Then, by the Hölder inequality, we obtain

(4.4)
$$\begin{aligned}|V_n(u) - V_n(v)| &\leq \frac{1}{H} \sum_{i=0}^{n-1} \Bigg( \left| \int_{t_i}^{t_{i+1}} (u_s - v_s) \, dB_s \right| \\ &\qquad \times \left| \int_{t_i}^{t_{i+1}} u_s \, dB_s \right|^{1/H} + \left| \int_{t_i}^{t_{i+1}} v_s \, dB_s \right| \Bigg)^{1/H} \\ &\leq \frac{1}{H} |V_n(u-v)|^H (V_n(u)^{1-H} + V_n(v)^{1-H}).\end{aligned}$$

If $u_t$ is a step process, then, from the properties of the fBm (see [15]), we know that $V_n(u)$ converges almost surely to $c \int_0^t |u_s|^{1/H} \, ds$. We can write

$$\begin{aligned}\left| V_n\left(\frac{1}{B}\right) - c \int_0^t \frac{ds}{|B_s|^{1/H}} \right| &\leq \left| V_n\left(\frac{1}{B}\right) - V_n(u_m) \right| + \left| V_n(u_m) - c \int_0^t |u_m(s)|^{1/H} \, ds \right| \\ &\quad + c \left| \int_0^t |u_m(s)|^{1/H} \, ds - \int_0^t \frac{ds}{|B_s|^{1/H}} \right| \\ &= a_{n,m} + b_{n,m} + c_m,\end{aligned}$$

where

$$u_m(t) = \sum_{j=0}^{m-1} \frac{1}{B_{t_{j-1}}} \mathbf{1}_{(t_{j-1}, t_j]}(t).$$

Using (4.4), we can estimate the first summand as follows:

$$a_{n,m} \leq \frac{1}{H} \left| V_n\left(\frac{1}{B} - u_m\right) \right|^H \left( V_n\left(\frac{1}{B}\right)^{1-H} + V_n(u_m)^{1-H} \right)$$

$$\leq \frac{1}{H}\left|\mathrm{Var}_{1/H}\left(\int_0^{\cdot}\left(\frac{1}{B}-u_m\right)dB_s\right)\right|^H$$
$$\times\left(\mathrm{Var}_{1/H}\left(\int_0^{\cdot}\frac{1}{B}dB_s\right)^{1-H}+\mathrm{Var}_{1/H}\left(\int_0^{\cdot}u_m\,dB_s\right)^{1-H}\right),$$

where $\mathrm{Var}_{1/H}$ denotes the total variation of order $\frac{1}{H}$ on the interval $[0,t]$. Then Love–Young's inequality (see [21]) implies that $\lim_m \sup_n a_{n,m}=0$. From the previous results, we know that, for any fixed $m$, $\lim_n b_{n,m}=0$, and clearly, $\lim_m c_m=0$. This completes the proof. □

We can now turn to the asymptotics of the skew-product decomposition. The main ingredient of what follows is the ergodic theorem.

PROPOSITION 16. *We have*

$$\frac{1}{\log t}\int_0^t \frac{\delta B_s}{B_s}\xrightarrow[t\to+\infty]{a.s.} H$$

*and*

$$\frac{1}{\log t}\int_0^t \frac{ds}{|B_s|^{1/H}}\xrightarrow[t\to+\infty]{a.s.} \frac{\Gamma(1-1/(2H))}{2^{1/(2H)}}.$$

PROOF. We first write

$$B_t = z_0 + \beta_t,$$

where $\beta$ is an fBm started at 0.

The key point is that

$$\frac{d\beta_s}{\beta_s} \quad \text{and} \quad \frac{ds}{|\beta_s|^{1/H}}$$

are invariant by scaling and this suggests to use Birkhoff–Khintchine's ergodic theorem.

Let $k > 0$, since the law of the standard fBm is invariant under the transformation $T_k\colon \omega \to \frac{\omega(k\cdot)}{k^H}$, by Birkhoff–Khintchine's theorem and Lemmas 12 and 14, we have

$$\frac{1}{N}\sum_{n=0}^{N-1}\int_1^k \frac{d\beta_{k^n s}/k^{nH}}{\beta_{k^n s}/k^{nH}}\xrightarrow[N\to+\infty]{a.s.} Y,$$

where $Y$ is constant. Hence,

$$\frac{1}{N}\sum_{n=0}^{N-1}\int_{k^n}^{k^{n+1}}\frac{d\beta_s}{\beta_s}\xrightarrow[N\to+\infty]{a.s.} Y.$$

We conclude that
$$\frac{1}{N}\int_1^{k^N} \frac{d\beta_s}{\beta_s} \xrightarrow[N\to+\infty]{\text{a.s.}} Y.$$

Since $k$ was arbitrary,
$$\frac{1}{\log t}\int_1^t \frac{d\beta_s}{\beta_s} \xrightarrow[t\to+\infty]{\text{a.s.}} Y.$$

Now, because $B$ is transient,
$$\frac{1}{\log t}\int_1^t \frac{d\beta_s}{\beta_s} \sim \frac{1}{\log t}\int_0^t \frac{dB_s}{B_s}.$$

This implies
$$\frac{1}{\log t}\int_0^t \frac{dB_s}{B_s} \xrightarrow[t\to+\infty]{\text{a.s.}} Y.$$

Let us now show that
$$\operatorname{Re}(Y) = H.$$

For this, we must check that
$$\mathbb{E}\left(\frac{\log |B_t|}{\log t}\right) \xrightarrow[t\to+\infty]{} H,$$
which is a direct consequence of the scaling property for $\beta$. We deduce that
$$\operatorname{Re} Y = H.$$

Now, we note that
$$\operatorname{Im} Y = 0$$
because of the rotational invariance of the fBm.

Since the second limit is obtained exactly by the same way, we omit the details of its proof. Let us just mention that it is based on the following computation:
$$\lim_{t\to+\infty} \mathbb{E}\left(\frac{1}{\log t}\int_0^t \frac{ds}{|B_s|^{1/H}}\right) = \mathbb{E}\left(\frac{1}{|N|^{1/H}}\right),$$
where $N$ denotes a random variable with standard normal distribution on $\mathbb{R}^2$. □

REMARK 17. Unlike what happens in the case of the standard Brownian motion, the process $(\int_0^t \frac{\delta B_s}{B_s})_{t\geq 0}$ is not a complex fBm $\beta$ time-changed with the clock
$$\int_0^t \frac{ds}{|B_s|^{1/H}}.$$

Indeed, otherwise, we would have

$$\frac{1}{\log t} \int_0^t \frac{\delta B_s}{B_s} = \frac{\beta_{\int_0^t ds/|B_s|^{1/H}}}{\log t} \xrightarrow[t \to +\infty]{\text{a.s.}} 0$$

because

$$\frac{\beta_t}{t} \xrightarrow[t \to +\infty]{\text{a.s.}} 0.$$

REMARK 18. It is also possible to apply the ergodic theorem to obtain, more generally, jointly with Proposition 11, the following limit theorem:

$$\frac{1}{\log t} \int_1^t \frac{1}{|\beta_s|^{1/H}} f\left(\frac{\beta_s}{|\beta_s|}\right) ds \xrightarrow[t \to +\infty]{\text{a.s.}} \frac{\Gamma(1 - 1/(2H))}{2^{1/(2H)}} \int_{\mathbb{S}} f \, d\sigma,$$

where $f$ is an integrable function on the unit circle $\mathbb{S}$ endowed with its normalized Haar measure $\sigma$ and $\beta$ is a two-dimensional fBm issued from 0.

As a direct corollary of this, and by time inversion we, deduce the following corollary.

COROLLARY 19. *Let $(B_t)_{t \geq 0}$ be a complex fBm with Hurst parameter $H > \frac{1}{2}$ and started at 0. We have, for any $z \in \mathbb{C}$,*

$$\frac{\log|B_t + z|}{\log t} \xrightarrow[t \to +\infty]{\text{a.s.}} H \quad and \quad \frac{\log|B_t + zt^{2H}|}{\log t} \xrightarrow[t \to 0^+]{\text{a.s.}} H.$$

PROOF. The first part of the proposition stems immediately from the skew-product decomposition, indeed,

$$\frac{\log|B_t|}{\log t} = \frac{1}{\log t} \operatorname{Re}\left(\int_0^t \frac{\delta B_s}{B_s}\right).$$

For the second part, we use a time inversion argument. Indeed, we can write

$$B_t = t^{2H} \beta_{1/t},$$

where $\beta$ is a planar fractional Brownian motion issued from 0, and then we apply the first part of the proposition to $\beta$. □

Observe that, from the law of iterated logarithm (see, e.g., [9] and [10]), it is known that, for any $H \in (0, 1)$, the following inequality holds:

$$|B_t| \leq A(1 + t^H)\sqrt{\ln\ln(3 + t)}.$$

Before we conclude this section, we give a small extension of Proposition 16.

PROPOSITION 20. *Let $(B_t)_{t\geq 0}$ be a complex fBm with Hurst parameter $H > \frac{1}{2}$ and started at $z_0 \neq 0$. Let $f : \mathbb{C} \to \mathbb{C}$ be a holomorphic function on the whole plane such that $f(0) = 0$. Then*

$$\frac{1}{\log t} \int_0^t f\left(\frac{1}{B_s}\right) \delta B_s \xrightarrow[t \to +\infty]{a.s.} H f'(0).$$

PROOF. Indeed, we can write for $z \in \mathbb{C}^*$ and some holomorphic function $g : \mathbb{C} \to \mathbb{C}$,

$$f\left(\frac{1}{z}\right) = \frac{f'(0)}{z} + \frac{1}{z^2} g'\left(\frac{1}{z}\right).$$

It stems from this that

$$\frac{1}{\log t} \int_0^t f\left(\frac{1}{B_s}\right) \delta B_s = \frac{f'(0)}{\log t} \int_0^t \frac{\delta B_s}{B_s} + \frac{1}{\log t} \int_0^t \frac{1}{B_s^2} g'\left(\frac{1}{B_s}\right) \delta B_s.$$

Now, thanks to Itô's formula,

$$\frac{1}{\log t} \int_0^t \frac{1}{B_s^2} g'\left(\frac{1}{B_s}\right) \delta B_s = \frac{1}{\log t}\left[g\left(\frac{1}{z_0}\right) - g\left(\frac{1}{B_t}\right)\right],$$

which allows to conclude according to Proposition 16. □

**5. Study of the windings.** In this section we shall again assume that $z_0 \neq 0$. We denote

$$\theta_t = \int_0^t \frac{B_s^2 \, dB_s^1 - B_s^1 \, dB_s^2}{|B_s|^2}.$$

From the skew-product decomposition of the two-dimensional fractional Brownian motion, the process $\theta$ characterizes the windings of $B$. Though the complete study of the asymptotics of $\theta$ is not yet well understood, we present some results. First, it is easy to study the windings of $B$ up to $2\pi$.

PROPOSITION 21. *We have*

$$\mathbb{E}(e^{in\theta_t}) \underset{n \to +\infty}{\sim} \frac{1}{t^{nH}}.$$

PROOF. From the skew-product decomposition, we get

$$\frac{B_t}{|B_t|} = \exp(i\theta_t), \qquad t \geq 0.$$

Hence, for $n \in \mathbb{N}$,

$$\mathbb{E}(\exp(in\theta_t)) = \mathbb{E}\left(\frac{(t^H N + z_0)^n}{|t^H N + z_0|^n}\right),$$

where $N$ is a two-dimensional standard normal law. The result follows then from a straightforward computation on the two-dimensional standard normal law. □

We conclude now the paper with a limit theorem for the functional of the two-dimensional fractional Brownian motion

$$(5.1) \qquad Z_t = \int_1^t \frac{B_s^2 \, dB_s^1 - B_s^1 \, dB_s^2}{s^{2H}},$$

which *looks* like the windings. Formula (2.3) allows us to compute the variance of this process:

$$\mathbb{E}(Z_t^2) = 2\alpha_H \int_1^t \int_1^t (rs)^{-2H} R(s,r)|r-s|^{2H-2} \, dr \, ds$$

$$(5.2) \qquad - 2\alpha_H^2 \int_{[1,t]^4} (sr)^{-2H} \mathbf{1}_{\{\theta \le r, \alpha \le s\}} |r-\alpha|^{2H-2} |\theta-s|^{2H-2} \, dr \, ds \, d\theta \, d\alpha$$

$$\le C_H \log t.$$

PROPOSITION 22. *Let $Z_t$ be defined by* (5.1). *Then $\frac{Z_t}{\sqrt{\log t}}$ converges in law, as $t$ tends to infinity, to the normal distribution $N(0, \sigma^2)$, where*

$$\sigma^2 = 4H(2H-1) \int_0^1 \left( y^{-2H} R(y,1) - H \log y^{-1} - H\beta\left(\frac{y}{1-y}\right) \right)(1-y)^{2H-2} \, dy$$

*and*

$$\beta(y) = \int_0^1 (1-x)^{2H-1}(x+y)^{-2H} dx.$$

PROOF. Define, for $\lambda \in \mathbb{R}$ and $t > 1$,

$$X_t = \exp\left(i\lambda \frac{Z_t}{\sqrt{\log t}}\right),$$

and set $\varphi(\lambda, t) = \mathbb{E}(X_t)$. Using the change of variable formula for the Riemann–Stieltjes integral we can write, for $1 < t_0 \le t$,

$$X_t = X_{t_0} + \int_{t_0}^t \frac{i\lambda}{\sqrt{\log u}} X_u u^{-2H} (B_u^2 \, dB_u^1 - B_u^1 \, dB_u^2)$$

$$(5.3) \qquad - \int_{t_0}^t \frac{i\lambda}{2u(\log u)^{3/2}} X_u Z_u \, du.$$

Using (2.4), we can transform this pathwise integral into a divergence integral, plus a complementary term:

$$\int_{t_0}^t \frac{i\lambda}{\sqrt{\log u}} X_u u^{-2H} (B_u^2 \, dB_u^1 - B_u^1 \, dB_u^2)$$

$$(5.4) \quad = \int_{t_0}^{t} \frac{i\lambda}{\sqrt{\log u}} X_u u^{-2H} (B_u^2 \, \delta B_u^1 - B_u^1 \, \delta B_u^2)$$

$$+ \alpha_H i\lambda \int_{t_0}^{t} \frac{u^{-2H}}{\sqrt{\log u}} \int_0^u (B_u^2 D_r^1 X_u - B_s^1 D_r^2 X_u)(u-r)^{2H-2} \, dr \, du.$$

Substituting (5.4) into (5.3) and taking expectations, we obtain

$$\varphi(\lambda, t) = \varphi(\lambda, t_0)$$

$$- \alpha_H \lambda^2 \int_{t_0}^{t} \frac{u^{-2H}}{\log u} \mathbb{E} \int_0^u X_u (B_u^2 D_r^1 Z_u - B_s^1 D_r^2 Z_u)(u-r)^{2H-2} \, dr \, du$$

$$- \int_{t_0}^{t} \frac{\lambda}{2u \log u} \frac{\partial \varphi}{\partial \lambda}(\lambda, u) \, du.$$

Now we apply the duality relationship (2.2) and we get

$$\mathbb{E}(X_u B_u^2 D_r^1 Z_u) = \mathbb{E}(X_u \langle D_\cdot^2 D_r^1 Z_u, \mathbf{1}_{(0,u)} \rangle_\mathcal{H})$$

$$+ \mathbb{E}(\langle D_\cdot^2 X_u, \mathbf{1}_{(0,u)} \rangle_\mathcal{H} D_r^1 Z_u)$$

$$= \mathbb{E}(X_u) \langle D_\cdot^2 D_r^1 Z_u, \mathbf{1}_{(0,u)} \rangle_\mathcal{H}$$

$$+ \frac{i\lambda}{\sqrt{\log u}} \mathbb{E}(X_u \langle D_\cdot^2 Z_u, \mathbf{1}_{(0,u)} \rangle_\mathcal{H} D_r^1 Z_u)$$

and we obtain

$$\varphi(\lambda, t) = \varphi(\lambda, t_0)$$

$$- \alpha_H \lambda^2 \int_{t_0}^{t} \frac{u^{-2H}}{\log u} \mathbb{E}(X_u)$$

$$\times \int_0^u (\langle D_\cdot^2 D_r^1 Z_u, \mathbf{1}_{(0,u)} \rangle_\mathcal{H}$$

$$- \langle D_\cdot^1 D_r^2 Z_u, \mathbf{1}_{(0,u)} \rangle_\mathcal{H})(u-r)^{2H-2} \, dr \, du$$

$$- \int_{t_0}^{t} \frac{\lambda}{2u \log u} \frac{\partial \varphi}{\partial \lambda}(\lambda, u) \, du + \psi(\lambda, t),$$

where

$$\psi(\lambda, t) = -i\lambda^3 \alpha_H \int_{t_0}^{t} \frac{u^{-2H}}{(\log u)^{3/2}}$$

$$\times \int_0^u \mathbb{E}[X_u(\langle D_\cdot^2 Z_u, \mathbf{1}_{(0,u)} \rangle_\mathcal{H} D_r^1 Z_u$$

$$- \langle D_\cdot^1 Z_u, \mathbf{1}_{(0,u)} \rangle_\mathcal{H} D_r^2 Z_u)](u-r)^{2H-2} \, dr \, du.$$

We have, for $r \leq u$,

$$D_r^1 Z_u = \frac{B_r^2}{r^{2H}} - \int_r^u \frac{\delta B_\theta^2}{\theta^{2H}},$$

$$D_r^2 Z_u = \int_r^u \frac{\delta B_\theta^1}{\theta^{2H}} - \frac{B_r^1}{r^{2H}}.$$

Hence,

$$\langle D_\cdot^1 D_r^2 Z_u, \mathbf{1}_{(0,u)}\rangle_\mathcal{H} = \langle \theta^{-2H} \mathbf{1}_{(r,u)} - r^{-2H}\mathbf{1}_{(0,r)}, \mathbf{1}_{(0,u)}\rangle_\mathcal{H}$$

$$= \alpha_H \int_0^u \int_r^u \theta^{-2H}|\sigma - \theta|^{2H-2}\, d\theta\, d\sigma - r^{-2H} R(r,u)$$

$$= H \int_r^u \theta^{-2H}(\theta^{2H-1} + (u-\theta)^{2H-1})\, d\theta - r^{-2H} R(r,u)$$

$$= H \log \frac{u}{r} + H\beta\left(\frac{r}{u-r}\right) - r^{-2H} R(r,u),$$

where

$$\beta(y) = \int_0^1 (1-x)^{2H-1}(x+y)^{-2H}\, dx,$$

and, similarly,

$$\langle D_\cdot^2 D_r^1 Z_u, \mathbf{1}_{(0,u)}\rangle_\mathcal{H} = \langle\, r^{-2H}\mathbf{1}_{(0,r)} - \theta^{-2H}\mathbf{1}_{(r,u)}, \mathbf{1}_{(0,u)}\rangle_\mathcal{H}$$

$$= r^{-2H} R(r,u) - H \log \frac{u}{r} - H\beta\left(\frac{r}{u-r}\right).$$

Hence,

$$\varphi(\lambda, t) = \varphi(\lambda, t_0)$$

$$- 2\alpha_H \lambda^2 \int_{t_0}^t \frac{u^{-2H}}{\log u} \mathbb{E}(X_u)$$

$$\times \int_0^u \left(r^{-2H} R(r,u) - H\log\frac{u}{r} - H\beta\left(\frac{u}{u-r}\right)\right)$$

(5.5)
$$\times (u-r)^{2H-2}\, dr\, du$$

$$- \int_{t_0}^t \frac{\lambda}{2u \log u}\frac{\partial \varphi}{\partial \lambda}(\lambda, u)\, du + \psi(\lambda, t)$$

$$= \varphi(\lambda, t_0) - 2\alpha_H \lambda^2 \int_{t_0}^t \frac{1}{u \log u}\varphi(\lambda, u)\rho(u)\, du$$

$$- \int_{t_0}^t \frac{\lambda}{2u \log u}\frac{\partial \varphi}{\partial \lambda}(\lambda, u)\, du + \psi(\lambda, t),$$

where
$$\rho(z) = \int_{1/z}^{1} \left( y^{-2H} R(y,1) - H \log y^{-1} - H\beta\left(\frac{y}{1-y}\right) \right)(1-y)^{2H-2} \, dy.$$

We claim that $\lim_{t\to\infty} \psi(\lambda, t) = 0$. In order to show this property, let us write

$$\mathbb{E}(X_u \langle D_\cdot^2 Z_u, \mathbf{1}_{(0,u)} \rangle_\mathcal{H} D_r^1 Z_u)$$
$$= \mathbb{E}\left( X_u \langle D_\cdot^2 Z_u, \mathbf{1}_{(0,u)} \rangle_\mathcal{H} \left( \frac{B_r^2}{r^{2H}} - \int_r^u \frac{\delta B_\theta^2}{\theta^{2H}} \right) \right)$$
$$= \mathbb{E}(\langle D_\cdot^2 Z_u, \mathbf{1}_{(0,u)} \rangle_\mathcal{H} (r^{-2H} \langle D_\cdot^2 X_u, \mathbf{1}_{(0,r)} \rangle_\mathcal{H} - \langle D_\cdot^2 X_u, \theta^{-2H} \mathbf{1}_{(r,u)} \rangle_\mathcal{H}))$$
$$= \frac{i\lambda}{\sqrt{\log u}} \mathbb{E}(X_u \langle D_\cdot^2 Z_u, \mathbf{1}_{(0,u)} \rangle_\mathcal{H} (r^{-2H} \langle D_\cdot^2 Z_u, \mathbf{1}_{(0,r)} \rangle_\mathcal{H}$$
$$- \langle D_\cdot^2 Z_u, \theta^{-2H} \mathbf{1}_{(r,u)} \rangle_\mathcal{H})).$$

Therefore,

$$\psi(\lambda, t) = -i\lambda^3 \alpha_H$$
$$\times \mathbb{E} \int_{t_0}^t \frac{u^{-2H}}{\log u^{3/2}} X_u$$
$$\times \int_0^u (\langle D_\cdot^2 Z_u, \mathbf{1}_{(0,u)} \rangle_\mathcal{H}$$
$$\times (r^{-2H} \langle D_\cdot^2 Z_u, \mathbf{1}_{(0,r)} \rangle_\mathcal{H} - \langle D_\cdot^2 Z_u, \theta^{-2H} \mathbf{1}_{(r,u)} \rangle_\mathcal{H})$$
$$+ \langle D_\cdot^1 Z_u, \mathbf{1}_{(0,u)} \rangle_\mathcal{H}$$
$$\times (r^{-2H} \langle D_\cdot^1 Z_u, \mathbf{1}_{(0,r)} \rangle_\mathcal{H} - \langle D_\cdot^1 Z_u, \theta^{-2H} \mathbf{1}_{(r,u)} \rangle_\mathcal{H}))$$
$$\times (u-r)^{2H-2} \, dr \, du.$$

We have

$$\langle D_\cdot^1 Z_u, \mathbf{1}_{(0,u)} \rangle_\mathcal{H} = \alpha_H \int_0^u \int_0^u \left( \frac{B_r^2}{r^{2H}} - \int_r^s \frac{\delta B_\theta^2}{\theta^{2H}} \right) |r-\theta|^{2H-2} \, dr \, d\theta.$$

Hence,

$$\|\langle D_\cdot^1 Z_u, \mathbf{1}_{(0,u)} \rangle_\mathcal{H}\|_2 \leq \alpha_H \int_0^u \int_0^u \left\| \frac{B_r^2}{r^{2H}} - \int_r^s \frac{\delta B_\theta^2}{\theta^{2H}} \right\|_2 |r-\theta|^{2H-2} \, dr \, d\theta$$
$$\leq C \int_0^u \int_0^u r^{-H} |r-\theta|^{2H-2} \, dr \, d\theta$$
$$\leq C u^H,$$

because, applying the inequality proved in [7], we obtain, for some constant $c_H$,

$$\left\|\int_r^u \frac{\delta B_\theta^2}{\theta^{2H}}\right\|_2 \leq c_H \|\theta^{-2H} \mathbf{1}_{(r,u)}\|_{L^{1/H}(\mathbb{R}_+)} \leq c_H r^{-H}.$$

Similarly,

$$\langle D^1_\cdot Z_u, \theta^{-2H} \mathbf{1}_{(r,u)} \rangle_{\mathcal{H}} = \alpha_H \int_0^u \int_r^u \left(\frac{B_\sigma^2}{\sigma^{2H}} - \int_\sigma^s \frac{\delta B_\theta^2}{\theta^{2H}}\right) |\sigma - \theta|^{2H-2} \theta^{-2H} \, d\theta \, d\sigma$$

and

$$\|\langle D^1_\cdot Z_u, \theta^{-2H} \mathbf{1}_{(r,u)} \rangle_{\mathcal{H}}\|_2 \leq C \int_0^u \int_r^u \sigma^{-H} |\sigma - \theta|^{2H-2} \theta^{-2H} \, d\theta \, d\sigma$$

$$= C d r^{-H},$$

where

$$d = \int_0^\infty \int_1^\infty y^{-H} |y - x|^{2H-2} x^{-2H} \, dx \, dy < \infty.$$

Also,

$$\|\langle D^1_\cdot Z_u, \mathbf{1}_{(0,r)} \rangle_{\mathcal{H}}\|_2 \leq \alpha_H \int_0^u \int_0^r \left\|\int_\sigma^s \frac{\delta B_\theta^2}{\theta^{2H}} - \frac{B_\sigma^2}{\sigma^{2H}}\right\|_2 |\sigma - \theta|^{2H-2} \, d\theta \, d\sigma$$

$$\leq C \int_0^u \int_0^r \sigma^{-H} |\sigma - \theta|^{2H-2} \, d\theta \, d\sigma$$

$$\leq C r^H.$$

Finally,

(5.6)
$$|\psi(\lambda, t)| \leq C \int_{t_0}^t \frac{u^{-H}}{\log u^{3/2}} \int_0^u r^{-H} (u-r)^{2H-2} \, dr \, ds$$

$$\leq \frac{C}{\sqrt{\log t}},$$

where the constant $C$ depends on $t_0$, $H$ and $\lambda$. Similarly, we can derive the estimate

(5.7)
$$\left|\frac{\partial \psi}{\partial t}(\lambda, t)\right| \leq \frac{C}{t(\log t)^{3/2}}.$$

Equation (5.5) yields

$$\frac{\partial \varphi}{\partial t} = -2\alpha_H \lambda^2 \frac{\varphi \rho}{t \log t} - \frac{\lambda}{2t \log t} \frac{\partial \varphi}{\partial \lambda} + \frac{\partial \psi}{\partial t}.$$

Hence,

$$\varphi(\lambda, t) = \varphi(\lambda, t_0) \exp\left\{-\int_{t_0}^t \frac{2\alpha_H \lambda^2 \rho(s)}{s \log s} ds\right\}$$

(5.8)
$$+ \int_{t_0}^t \exp\left\{-\int_{t_0}^s \frac{2\alpha_H \lambda^2 \rho(u)}{u \log u} du\right\} \left(-\frac{\lambda}{2s \log s} \frac{\partial \varphi}{\partial \lambda} + \frac{\partial \psi}{\partial s}\right) ds.$$

The integral $\exp\{-\int_{t_0}^t \frac{2\alpha_H \lambda^2 \rho(s)}{s \log s} ds\}$ behaves as $\frac{1}{\log t}$ as $t$ tends to infinity. We have, from (5.2),

(5.9)
$$\left|\frac{\partial \varphi}{\partial \lambda}\right| \leq \frac{|\lambda|}{\sqrt{\log t}} |\mathbb{E}(Z_t X_t)| \leq |\lambda| \sqrt{C_H}.$$

Then, (5.9) and (5.7) imply that $\lim_{t \to \infty} \varphi(\lambda, t) = \varphi(\lambda, \infty)$ exists. By a similar argument, differentiating equation (5.8) with respect to $\lambda$, we can show that $\lim_{t \to \infty} \frac{\partial \varphi}{\partial \lambda} = \frac{\partial \varphi}{\partial \lambda}(\lambda, \infty)$ exists. We claim that

(5.10)
$$-2\alpha_H \lambda \rho(\infty) \varphi(\lambda, \infty) - \frac{1}{2} \frac{\partial \varphi}{\partial \lambda}(\lambda, \infty) = 0.$$

Otherwise, $\frac{\partial \varphi}{\partial t}$ behaves as $\frac{C}{t \log t}$ when $t$ tends to infinity, where $C$ is a constant different from zero, which is contradictory because we would have

$$\varphi(\lambda, t) \underset{t \to +\infty}{\sim} C \log \log t.$$

Equality (5.10) implies that

$$\varphi(\lambda, \infty) = e^{-2\alpha_H \lambda^2 \rho(\infty)},$$

and $\frac{Z_t}{\sqrt{\log t}}$ converges in distribution as $t$ tends to infinity to the normal distribution $N(0, 4\alpha_H \rho(\infty))$. □

Notice also that, as $H$ tends to $\frac{1}{2}$, the variance $\sigma^2$ in the preceding proposition converges to 2. This is reasonable, because if $B_t$ is a standard two-dimensional Brownian motion, $\frac{1}{\sqrt{\log t}} \int_1^t \frac{B_s^2 dB_s^1 - B_s^1 dB_s^2}{s}$ converges in distribution as $t$ tends to infinity to the normal law $N(0, 2)$. Indeed, the Girsanov theorem implies that

$$\mathbb{E}\left(\exp\left(\frac{i\lambda}{\sqrt{\log t}} \int_1^t \frac{B_s^2 dB_s^1 - B_s^1 dB_s^2}{s} + \frac{\lambda^2}{2 \log t} \int_1^t \frac{|B_s|^2}{s^2} ds\right)\right) = 1,$$

and, applying the ergodic theorem, we have that $\frac{1}{2 \log t} \int_1^t \frac{|B_s|^2}{s^2} ds$ converges almost surely to one. Nevertheless, we could expect that, in the case of the fBm which is transient, the functional

$$\int_1^t \frac{B_s^2 dB_s^1 - B_s^1 dB_s^2}{s^{2H}}$$

is *closer* to the windings of $B$ than it is the case for the planar Brownian motion for which it is well known, since [18], that they behave like $\frac{1}{2}\mathcal{C} \log t$, where $\mathcal{C}$ is a Cauchy law with parameter 1.

**Acknowledgments.** This work was initiated while the first author was visiting the IMUB from September to December 2002. He wishes to thank friendly D. Nualart for his hospitality, his enthusiasm and his financial support.

LABORATOIRE DE STATISTIQUE
ET PROBABILITÉS
UNIVERSITÉ PAUL SABATIER
118 ROUTE DE NARBONNE
31500 TOULOUSE
FRANCE
E-MAIL: fbaudoin@cict.fr

FACULTAT DE MATEMÀTIQUES
UNIVERSITÉ DE BARCELONA
GRAN VIA 585
08007 BARCELONA
SPAIN
E-MAIL: dnualart@ub.edu